\documentclass[a4paper, 
11pt
]{article}

\usepackage{subfig}
\usepackage{pgfplots}
\usepackage{pgfplotstable}
\usepackage{mathtools}
\usepackage{multicol}
\usepackage{comment}
\usepackage{booktabs}
\pgfplotsset{compat=1.5}
\usepackage{amssymb}
\usepackage{url}

\usepackage{pdfsync}
\usepackage{float}
\usepackage{tabularx}
\usepackage{enumerate}
\usepackage{array}
\usepackage{xspace}

\usepackage{authblk}
\usepackage{tikz}
\usepackage{tikzsymbols}
\usetikzlibrary{calc,trees,positioning,arrows,chains,shapes.geometric,%
    decorations.pathreplacing,decorations.pathmorphing,shapes,%
    matrix,shapes.symbols, decorations.markings, patterns,fit}

\newcommand{\mupar}{\ensuremath{\boldsymbol{\mu}}}
\newcommand{\etapar}{\ensuremath{\boldsymbol{\eta}}}

\begin{document}

\title{An integrated data-driven computational pipeline with model order reduction for
  industrial and applied mathematics}

\author[]{Marco~Tezzele\footnote{marco.tezzele@sissa.it}}
\author[]{Nicola~Demo\footnote{nicola.demo@sissa.it}}
\author[]{Andrea~Mola\footnote{andrea.mola@sissa.it}}
\author[]{Gianluigi~Rozza\footnote{gianluigi.rozza@sissa.it}}

\affil[]{Mathematics Area, mathLab, SISSA, International School of Advanced Studies, via Bonomea 265, I-34136 Trieste, Italy}

\maketitle

\begin{abstract}
In this work we present an integrated computational pipeline involving several
    model order reduction techniques for industrial and applied mathematics, as
    emerging technology for product and/or process design procedures. Its
    data-driven nature and its modularity allow an easy integration into
    existing pipelines. We describe a complete optimization framework with
    automated geometrical parameterization, reduction of the dimension of the
    parameter space, and non-intrusive model order reduction such as dynamic
    mode decomposition and proper orthogonal decomposition with interpolation.
    Moreover several industrial examples are illustrated.
\end{abstract}

\section{Introduction}
\label{sec:intro}
A very common problem in the optimization of the design of
industrial artifacts is that of finding the shape that minimizes some
quantity of interest representing the expected performance. From a
mathematical point of view such a problem translates into an
optimization problem in which a suitable algorithm makes several
queries to a simulation solver allowing for an evaluation of each
sample in the design space. This leads to the identification of the optimal
solution, possibly subjected to a set of prescribed constraints.

The experience acquired through several industrial projects suggested us 
that for such pipeline to operate in a robust way in a manufacturing
environment, several aspects have to be integrated and developed so as
to deal with both complex geometries and solution fields. With the
concept of digital twin becoming widespread nowadays, we have to be able to
pass automatically from industrial CAD geometries to fluid dynamics and
structural simulations which allow for virtual performance
evaluation. All the steps in the procedure that, moving from a CAD
geometry, leads to an optimized shape need to be carefully devised and
integrated. First, one has to process the industrial geometry at hand
through a suitable shape parameterization strategy which identifies
the parameter space describing all possible designs to be
investigated. After the generation of a suitable space including all the parameters that
satisfy all the structural, geometrical, and regulatory constraints
prescribed by the design engineers according to the stakeholders
needs, a proper sampling of such space is used to set up a campaign of
CFD and structural simulations resulting in the performance evaluation
of each shape tested. Usually, a single full order industrial simulation takes days
or even weeks, so it is crucial to develop reduced order models (ROMs) so as
to speed up the full optimization cycle and make it compatible
with the production needs. The computational time of a single sample point
evaluation is reduced through different ROM techniques such as dynamic
mode decomposition (DMD), and proper orthogonal decomposition
(POD), both based on singular value
decomposition~\cite{benner2017model,benner2017morepas,salmoiraghi2016advances,rozza2018advances}. In
the case of DMD, ROM is used to reconstruct and predict the
spatiotemporal dynamics of a high fidelity simulation such that its
evolution can be completed at a faster rate. Instead POD-Galerkin or
POD with interpolation exploit data on previous simulations, properly
stored, to provide accurate surrogate solutions corresponding to
untested sample points in the parameter space. In such a way the
computational cost of an online optimization cycle can be dropped to
hours or even minutes.

In a further post-processing phase, we also apply a reduction of the
parameter space exploiting the active subspaces
property~\cite{constantine2015active,tezzele2018dimension,tezzele2018combined}. Such an 
analysis allows for the detection of possible redundancies in the
chosen parameters, suggesting linear combinations of the original
parameters which dominate the system response.

This work aims at presenting and discussing a series of best practice approaches
for the application of each of the aforementioned techniques within an industrial
optimization framework. Such approaches are the result of the constant involvement
of mathLab laboratory of SISSA\footnote{\url{www.mathlab.sissa.it}} in
industrial projects joining the research efforts of both manufacturing
companies and academic institutions.  After a brief overview of the overall
problem and goals, the contribution is arranged so as to present each of the
described industrial numerical pipeline steps in a complete and detailed
fashion. We first consider the geometrical treatment of the industrial
artifacts shape, then we suggest possible parameter space reduction strategies.
We then provide details on data-driven model order reduction methodologies, and
finally present a brief summary of numerical results obtained in some
applications carried out in the framework of industrial projects in which
mathLab is involved.

\section{From digital twin to real-time analysis}
Nowadays the digital twin is a concept well spread among all the
engineering companies and communities. With the exact digital representation of a
physical system, it allows, for example, to perform structural and
fluid dynamics simulations, to make sensitivity analysis with respect
to the parameters, and to optimize the design. The increased number of
devices for real-time data acquisition of the physical system makes the
digital twin paradigm obsolete, due to its intrinsic static nature. 

We are moving every day into a more dynamical representation of the entire
system that takes into account more and more real-time data from different
sources. This new paradigm, thanks to data-driven models, uncertainty
quantification, machine learning algorithms, artificial intelligence and better
integration of all the singular computational modules, will provide new
capabilities in terms of discovery of hidden correlations, fault detection,
predictive maintenance and design optimization. Its goal is to enable delivery
of better simulation and modeling results, and thus shorter the product
development, the so-called time-to-market, and reduce the product maintenance
and potential downtime.

These new needs from the industrial point of view, lead to new
mathematical methods and new interdisciplinary pipelines for data acquisition, model
order reduction, data elaboration, as well as optimization
cycles~\cite{abisset2018model,demo2018shape,ghnatios2012proper}.

\begin{figure}[htbp]
\centering
\begin{tikzpicture}[%
    scale=0.2,%
    >=triangle 60,              
    every join/.style={norm},   
    ]
\tikzset{
  base/.style={scale=0.8, draw, on grid, text width=6.em, align=center, minimum height=4ex},
  proc/.style={base, rectangle, text width=8em},
  test/.style={base, diamond, aspect=2, text width=3em},
  term/.style={proc, rounded corners, text width=5em}, 
  coord/.style={coordinate, on grid, node distance=6mm and 25mm},
  nmark/.style={draw, cyan, circle, font={\sffamily\bfseries}},
  norm/.style={->, draw, lcnorm},
  free/.style={->, draw, lcfree},
  cong/.style={->, draw, lccong},
  it/.style={font={\small\itshape}}
}
\node [term] (p0) {CAD file};
\node [term, below of=p0, node  distance=1.6cm] (p0a)
{Structural\\and CFD\\team};
\node [term, below of=p0a, node distance=1.8cm] (p0b) {Data\\acquisition};
\node [term, below of=p0b, node distance=1.4cm] (p0c) {Regulations};
\node [base, label={\small \it input}, below right=1cm and 3.2cm of
p0] (p1) {Geometrical, \\structural, \\ physical \\parameters};
\node [base, right of=p1, node distance=3.2cm] (p2) {Complete\\design};
\node [base, above of=p2, node distance=1.cm, pattern=north west lines, pattern color=orange!30] (p2a) {Intrusive};
\node [base, below of=p2, node distance=1.2cm, pattern=north west lines, pattern color=blue!30] (p2b) {Non-intrusive};
\node [base, below of=p2b, node distance=1.cm] (p2bb) {Home-made\\software};
\node [base, label={\small \it output}, right of=p2, node
distance=3.2cm] (p3) {Validation\\ and control};
\node [term, right of=p3, node distance=3.cm] (p4) {Final design};
\node[draw,dotted,inner sep=8pt,fit=(p2) (p2a) (p2bb), label={\bf offline-online}] {};
        \node[draw,dashed,inner sep=8pt,fit=(p2) (p2a) (p2bb),
        label={below:{\it \small black-box}}] {};

\draw [->] (p1) -- (p2);
\draw [->] (p2) -- (p3);
\draw [->] (p3) -- (p4);
\draw [->] (p3) -- (41.6,-23);
        \draw [-] (41.6,-23) -- node[above]{\it \small design feedback} (8,-23);
\draw [>-*] (p0) -- (8.4,0);
\draw [>-*] (p0a) -- (8.4,-6.5);
\draw [>-*] (p0b) -- (8.4,-13.5);
\draw [>-*] (p0c) -- (8.4,-19);
\draw [-] (8, -23) -- (8,0);
\draw [->] (8, -5) -- (p1);
\draw[decorate,decoration={brace,amplitude=8pt,mirror}] (9, -23.5) -- (21.5, -23.5) node [midway,below, yshift=-10pt] {};
        \draw[decorate,decoration={brace,amplitude=8pt,mirror}] (22.5, -23.5) -- (40.6, -23.5) node [midway,below, yshift=-10pt] {};
\node[base, draw=none, text width=12em] (aa) at (15, -27) {Parameter space\\reduction};
\node[base, draw=none] (bb) at (31, -27) {Computational \\reduction};
\end{tikzpicture}
\caption{Complete optimization pipeline involving automatic interface
  with CAD files, experimental and numerical data acquisition,
  definition of the parameters constraints, and parameter space
  dimensionality reduction. Two possible reduced order methods are
  presented: an intrusive and a non-intrusive approach, allowing the creation
  of a complete simulation-based design optimization framework.}
\label{fig:schema}
\end{figure}
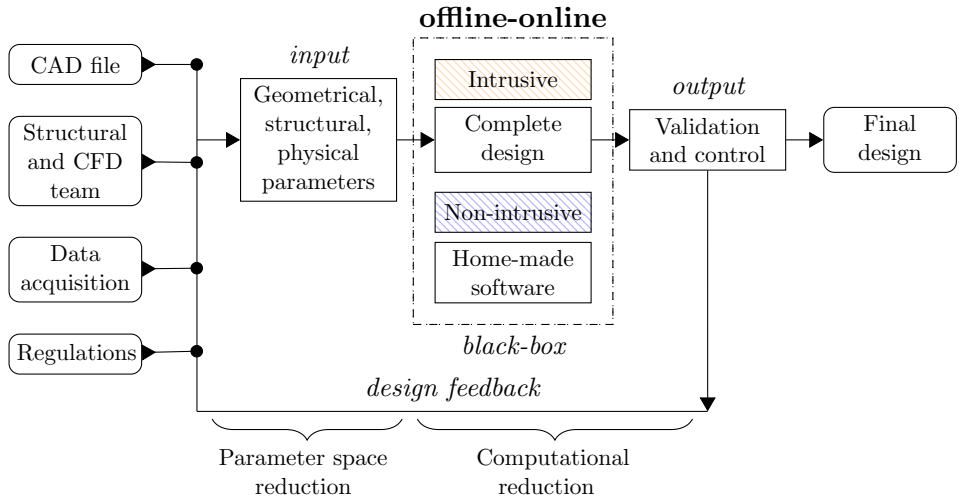

We can summarize a modern shape design optimization cycle with the diagram in
Figure~\ref{fig:schema}. Usually a CAD file describing the
geometry to be optimized is provided by the design team, then we have the
structural and CFD teams that provide physical constraints and
admissible range for the parameters variation, we have data coming
from the experiments on the scale model, and finally we have
regulations constraints from the national authorities. 

When the inputs of the simulation are set, the output is computed
through high fidelity solver and the optimization cycle is closed
by validation and control, the imposition of the constraints, and a
selection of a new set of parameters. If a complete simulation
takes several hours, even days, finding the optimal shape becomes impossible
due to the many evaluations of the parametric PDE needed by the
optimizer. Here, the model order reduction (MOR) approach allows fast
evaluations of the output of interest or derived quantities of the
output thus enabling to find the optimum shape in few hours and to
test different optimization algorithms. MOR techniques are very
versatile, enabling both intrusive and non-intrusive approach with
respect to the solver used. If it is a commercial black box it is
possible to use numerical methods that work using only precomputed
solution fields. While having access to the solver code allows also to
reduce the single operators of the PDEs.

After the continued validation and control the optimization cycle ends by
providing the final design in a suitable file format that can be
analysed by all the interested departments.

In the following we are going to present all the techniques and
integrated pipelines developed to accomplish such simulation based design
optimization.

\section{Advanced geometrical parametrization with automatic CAD files interface}
\label{sec:geom}
The first important step is related to the geometry of the shapes considered.
As previously mentioned, one of the aspects typically subjected to
optimization in the design of industrial artifacts is in fact their
geometry. Finding the shape by maximizing the performance of a certain
product or of one of its components is in fact a very common problem
in industry. From a mathematical standpoint, such class of problems
is obviously formalized as an optimization problem, which consists in
the identification of the point of a suitable parameter space that
maximizes the value of a prescribed performance parameter (or output
function). Although the mathematical algorithms carrying out such task
are commonly well assessed, the mathematical formalization of  the
problem requires that the shape modification can be recast as the
corresponding variation of a certain set of parameters. The latter
operation, which somehow translates the properties of the geometrical
shapes into a set of numbers handled by the optimization algorithms,
is usually referred to as \emph{shape parametrization}. In the most
common practice the shape parametrization is a rather delicate part of
the overall design process. In fact, as the shape to be studied can be
specified in several different file formats or analytical
descriptions, a unified approach for shape parametrization algorithms
has currently not been established. In the present section we will
describe and discuss the state of the art of shape parametrization
techniques, and present examples of their application to the geometry
of different industrial artifacts. 

\begin{figure}
\centering
\includegraphics[width=0.65\textwidth]{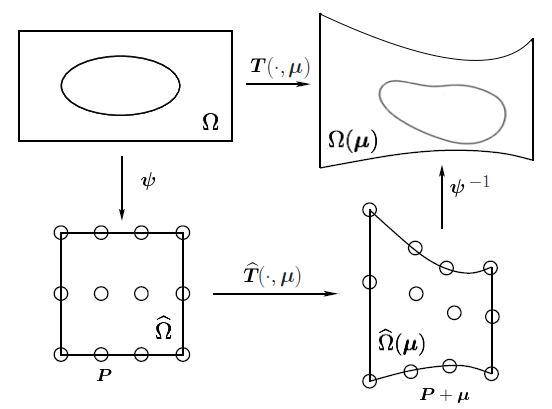}
\caption{Sketch of the three maps composing the FFD map
  construction. $\psi$ maps the physical space to the reference one,
  then $\hat{T}$ deforms the entire geometry according to the
  displacements of the lattice control points, and finally $\psi^{-1}$
maps back the reference domain to the original one.}
\label{fig:FFD_sketch}
\end{figure}

A first shape parametrization algorithm which has been devised so as
to be applied to arbitrarily shaped geometries, is the free form
deformation
(FFD)~\cite{sederbergparry1986,LassilaRozza2010,sieger2015shape}. FFD
consists basically in three different steps, as depicted in Figure~\ref{fig:FFD_sketch}. The first step is that of
mapping the physical domain into a reference one. In the second step,
a lattice of points is built in such reference domain, and some points
are moved to deform the lattice. Since the lattice points represent
the control points of a series of shape functions (typically Bernstein
polynomials), the displacement of such points is propagated to compute
the deformation of all the points of the reference domain. In the
third, final step, the deformed reference domain is mapped back into
the physical one, to obtain the resulting morphed geometry.

So, the displacements of one or more of the control points in the
lattice are the parameters of the FFD map, which is the composition of
the three maps described. As both the number of points in the lattice
and the number of control points displaced to generate the deformation
are flexible, the FFD map can be built with an arbitrary number of
parameters. Thus, a very useful feature of FFD, is that it allows for
the generation of global deformations even when few parameters are
considered. 

\begin{figure}[ht]
\centering
\subfloat[Original]{\includegraphics[width=0.48\textwidth]{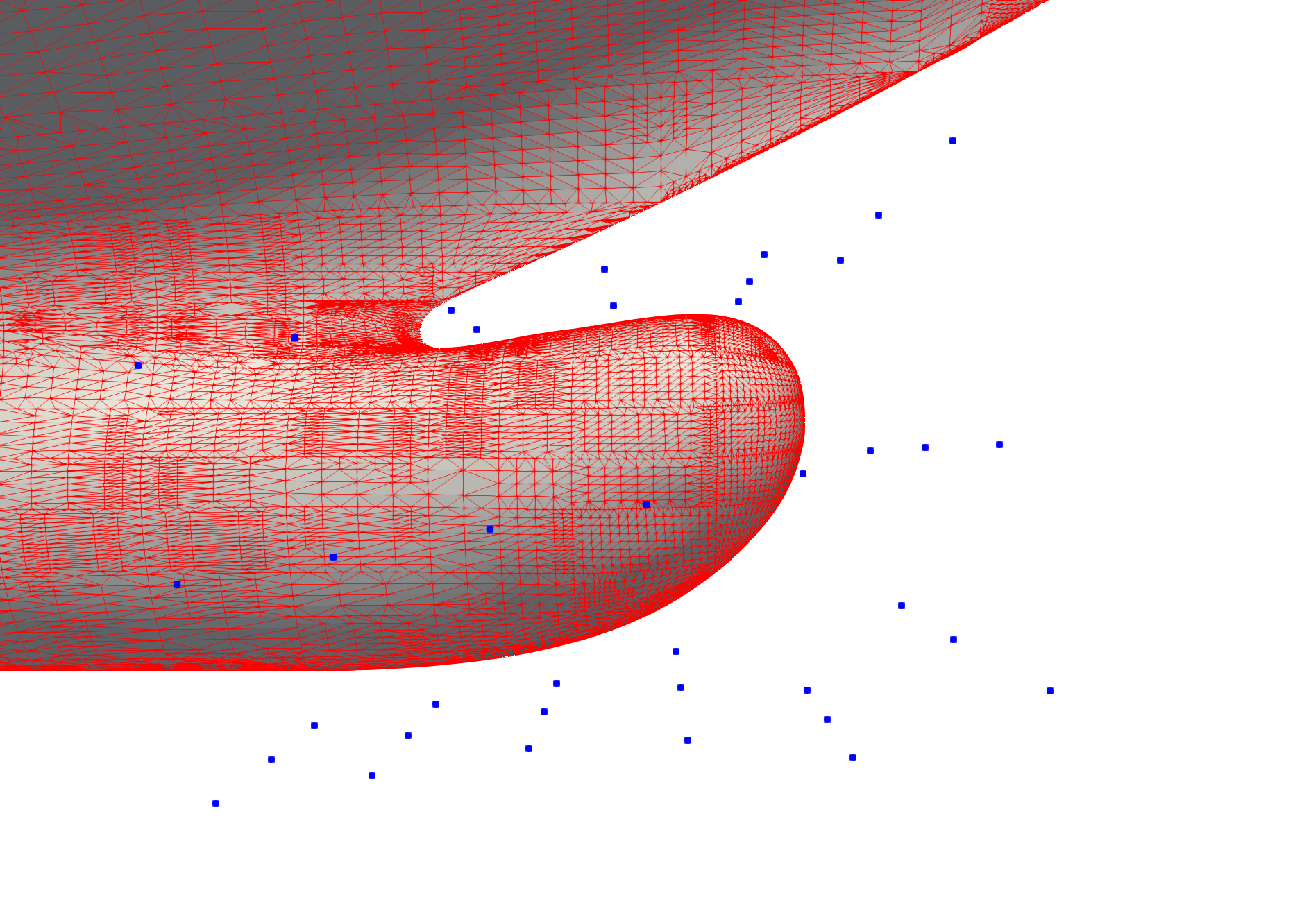}\label{subfig:ffd_hull_1}}
\hfill
\subfloat[Deformed]{\includegraphics[width=0.48\textwidth]{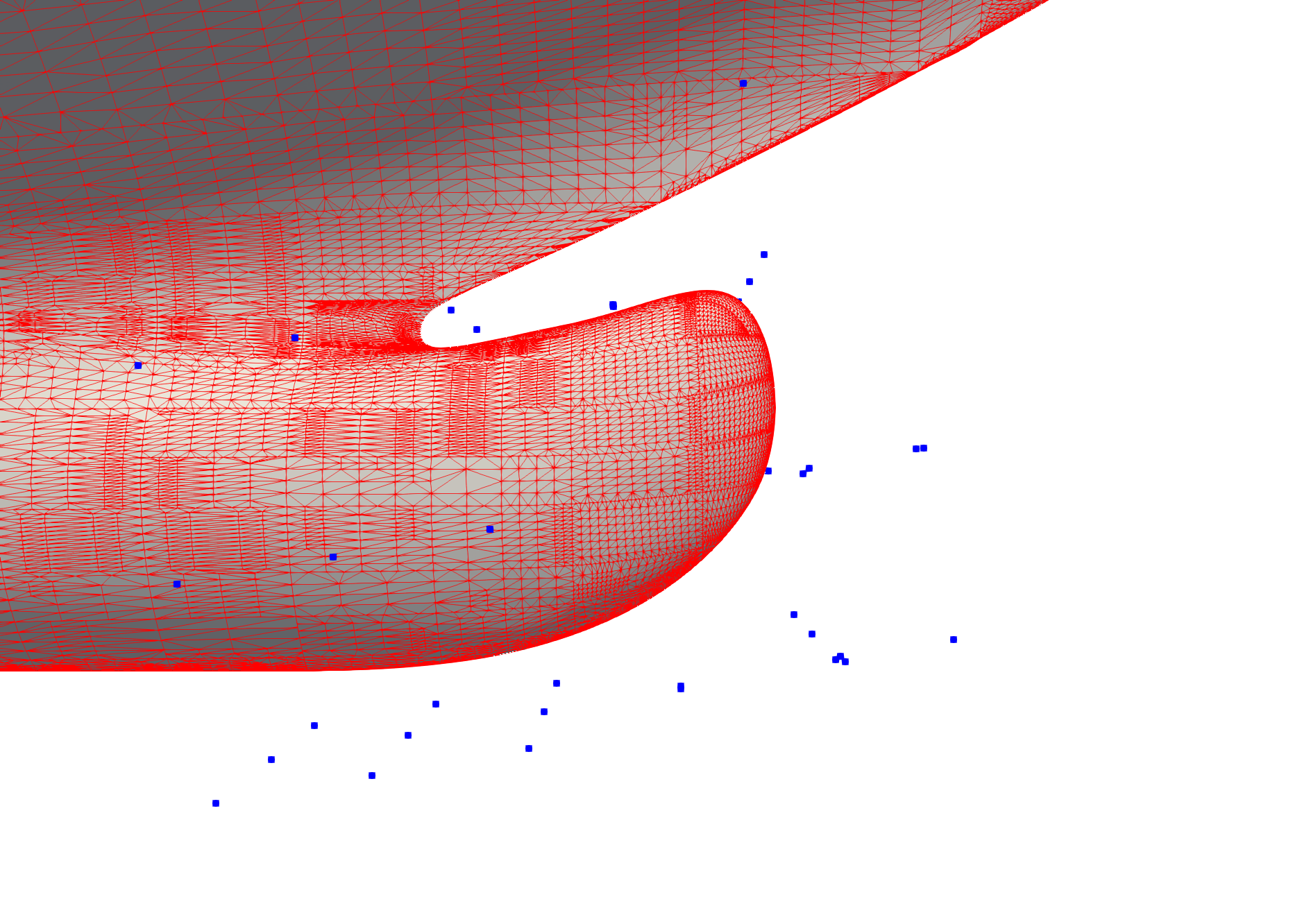}\label{subfig:ffd_hull_2}}
\caption{An example of the application of FFD to an STL triangulation. Plot \protect\subref{subfig:ffd_hull_1} shows the FFD lattice over
  one side of the bulbous bow of a ship. Plot \protect\subref{subfig:ffd_hull_2}
  depicts a deformed configuration of the same hull, along with the displaced FFD lattice control points.}
\label{fig:ffd_lattice}
\end{figure}

Since FFD is able to define a displacement law for each point of the
3D space contained in the control points lattice, it can be quite
naturally applied to all the geometries that are specified through
surface triangulations, surface grids or even volumetric grids. As a
first example, in Figure~\ref{fig:ffd_lattice} we present the
application of FFD to an STL triangulation, which is a very common
output format for CAD modeling tools. The shape deformation
illustrated in the picture has been carried out making use of the
PyGeM~\cite{pygem} Python package, which has been developed to be
directly interfaced with the industrial geometry files and to deform
them, so as to cut the communication times between the company
simulation team and design team. In this application, the STL
triangulation is imported and the coordinates of the nodes in the
triangulation are modified according to the FFD map generated through
the user specified lattice (included in the pictures).

\begin{figure}[ht]
\centering
\subfloat[Original]{\includegraphics[width=0.48\textwidth]{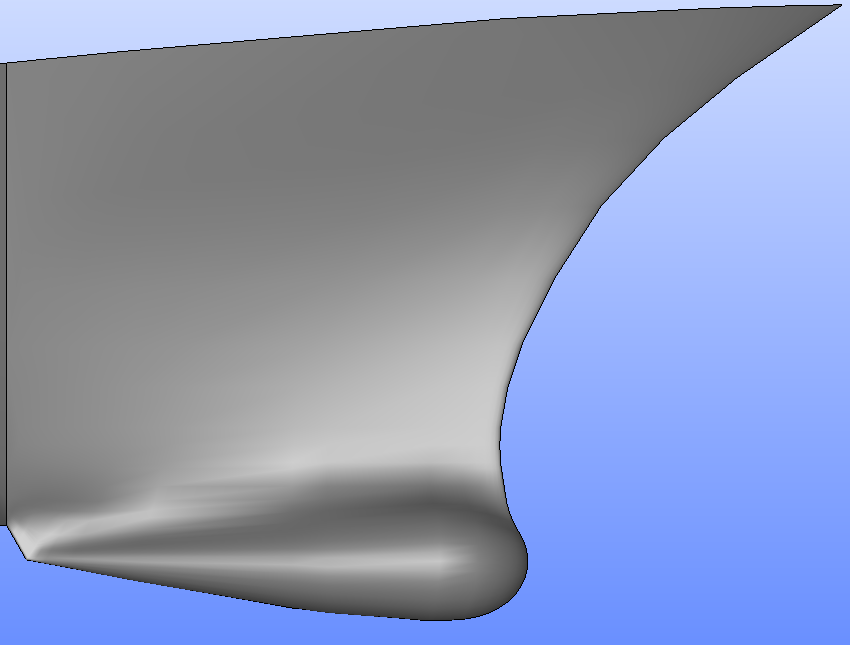}\label{subfig:ffd_dtmb_1}}
\hfill
\subfloat[Deformed]{\includegraphics[width=0.48\textwidth]{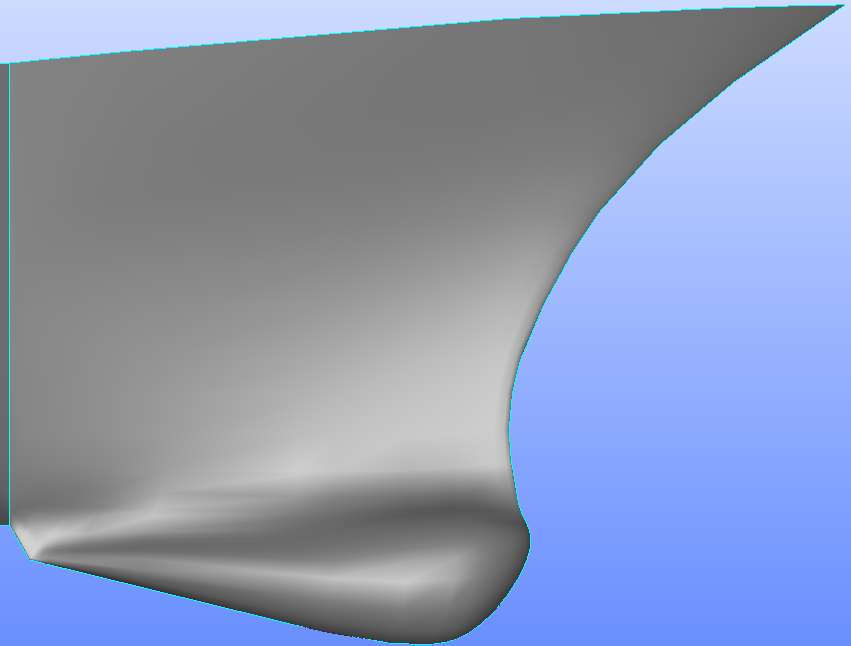}\label{subfig:ffd_dtmb_2}}
\caption{An example of the application of FFD to an IGES CAD geometry containing parametric surfaces. Plot \protect\subref{subfig:ffd_dtmb_1}
  shows the original shape of the DTMB-5415 US Navy Combatant hull bulb. Plot \protect\subref{subfig:ffd_dtmb_2}
  presents a deformed configuration of the same bulb.}
\label{fig:ffd_dtmb}
\end{figure}

The versatility of FFD can be further exploited to deform CAD surfaces
that are generated as the collection of parametric patches. In such
case, the desired deformation is obtained applying the specified FFD
map to the control points of the NURBS and B-Splines surfaces of each
patch composing the CAD model. The result in this case will also be a
deformed geometry which follows the FFD deformation law specified by
the user. Figure~\ref{fig:ffd_dtmb} presents a real world application
of FFD to CAD parametric surfaces. The renderings show the original
bulbous bow of the DTMB-5415 US Navy Combatant hull (which is a common
benchmark for the CFD simulations community), and one of its several
deformations generated to carry out the campaign of fluid dynamic
simulations discussed in~\cite{demo2018isope}. The PyGeM capabilities
allow for importing the CAD geometry (in IGES or STEP format), extract
and modify the control points of its surfaces and curves, and export
the result into new CAD files.

\begin{figure}[ht]
\centering
\includegraphics[width=1\textwidth]{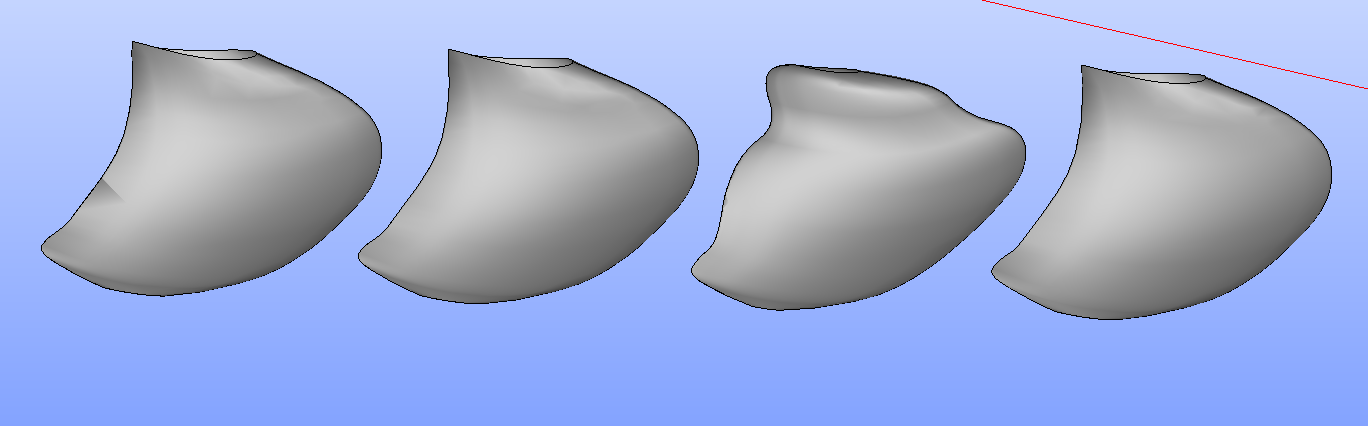}
\caption{An example of the deformation of an engineered propeller
  blade shape. The picture shows four deformed configurations of a
  PPTC benchmark propeller blade, in which the sectional properties of
  the blade have been kept untouched, while modifying the pitch, rake
  and skew distributions. The deformations are performed with the
  BladeX Python package~\cite{bladex}.}
\label{fig:ffd_pptc}
\end{figure}

Along with FFD, the PyGeM package implements other morphing techniques: the
deformation based on radial basis function (RBF)
interpolation~\cite{buhmann2003radial,morris2008cfd,manzoni2012model},
and the inverse distance weighting (IDW)
interpolation~\cite{shepard1968,witteveenbijl2009,forti2014efficient,BallarinDAmarioPerottoRozza2017}. Yet,
there are situations in which the shape to be deformed is already been
engineered in a specific way, and general purpose deformation
algorithms like FFD and the ones just mentioned would not preserve some
critical characteristics of the geometry. A rather striking example of
this is represented by the deformation of propeller blades illustrated in
Figure~\ref{fig:ffd_pptc}. The shape of a propeller blade is in fact
generated in a bottom-up fashion, first defining a set of sections
represented by airfoils, and then properly placing each section in the
three dimensional space. Since the aerodynamic feature of each airfoil
section are known to the engineers which have selected them, any
deformation that alters the sectional shape of the blade will not lead
to an acceptable geometry. Thus, for such highly engineered shapes
\emph{ad hoc} solutions have to be generated. In most cases such
solutions exploit the very algorithms used by the engineers to generate the
artifacts in the first place, introducing parameters in one or more
points of the generation procedure so as to create a class of
deformed shapes. Also in this case, once the shape has been properly
deformed (or re-generated), it has to be saved in the proper CAD
format (IGES, STEP or STL) to be handled by the CFD or design team.

\begin{figure}[h!]
\centering
\includegraphics[width=0.75\textwidth]{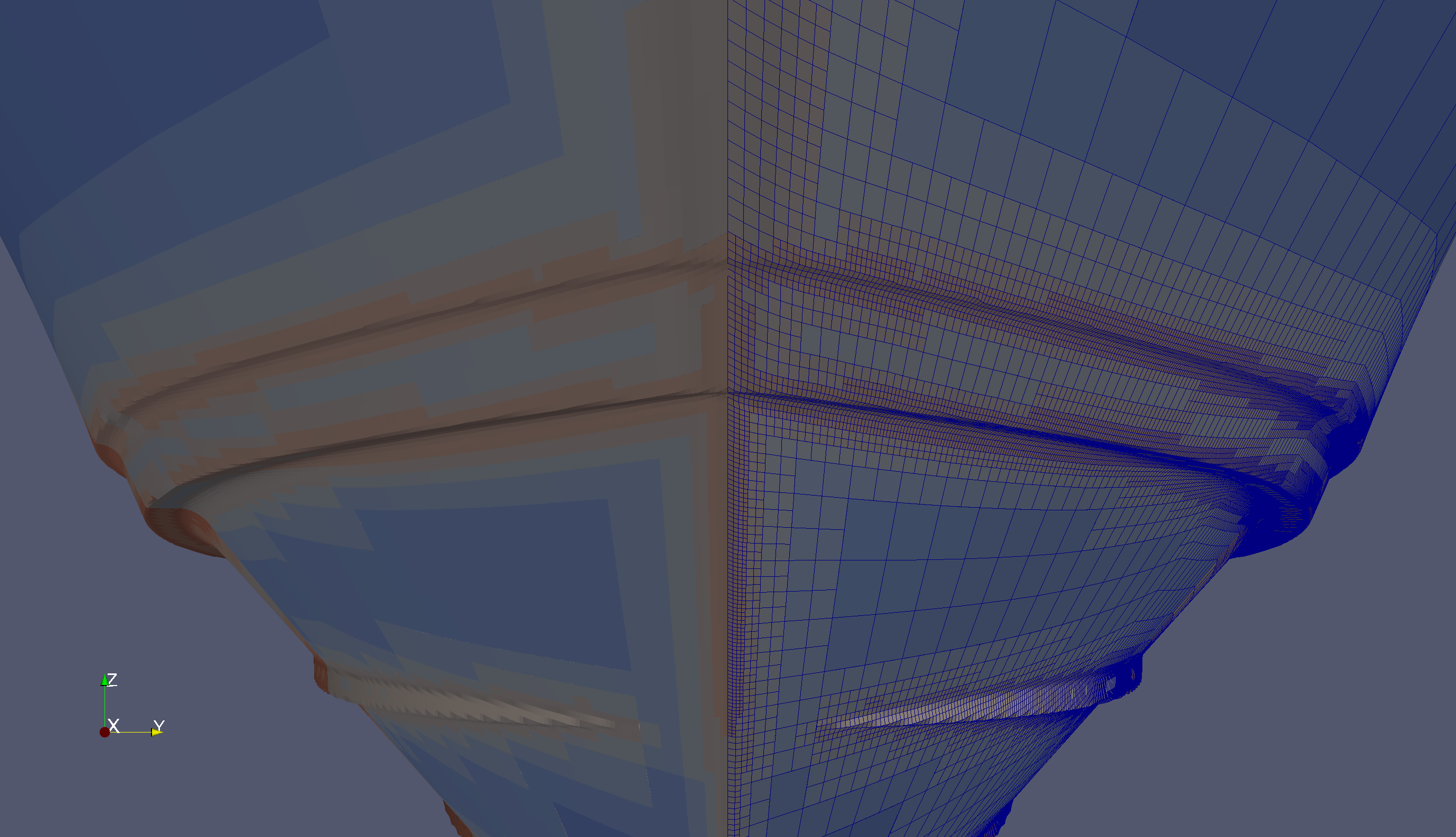}
\caption{Water tight quadrilateral mesh generated on the water non tight IGES surface of a planing yacht hull.}
\label{fig:hull_front}
\end{figure}

As mentioned, in the CAD data structures surfaces can be described both by means
of a triangulation and as a collection of parametric patches. Yet,
most CAD modeling tools
used by engineers to generate the virtual model of any object designed, operate
using NURBS and B-Splines curves and surfaces. On the other hand, several tools
for CFD or structural analysis only handle geometries specified through triangulations.
For such reason, a series of algorithms that generate triangulations on arbitrary
parametric surfaces has been implemented in the last years. Among others, the ability to
produce closed --- or \emph{water tight} ---  triangulations on surfaces composed of possibly
unconnected faces is a crucial feature for both CFD and structural engineering applications. 
In fact, the neighboring parametric patches composing a CAD model are generally only connected
up to a specified tolerance. Thus, it is not infrequent that ideally continuous surfaces present
small gaps and overlaps between each patch composing them. To avoid the problem, most CAD modeling tools
retain several logical information to complement the geometric data and indicate which patches should
be considered as neighbors. Yet, converting to vendor-neutral file formats such as IGES of STEP that
allow the digital exchange of information among CAD systems will cause in most case the loss of topological
information on neighboring patches. This is often a problem in the numerical analysis community, in which
geometries are in most cases obtained by third parties, and in which water tight surfaces are in needed to
define (and confine) the three dimensional domains considered in the simulations. So, a possible strategy
to avoid a surface triangulation that depends on the local patches parametrization, and suffers from
parametrization jumps and gaps, is to create new nodes not in the surface parametric space, but in
the physical three dimensional space. Since such new nodes will not be initially located onto the
CAD surface, a series of surface and curves projectors are used to make sure that the new grid points
are properly placed onto the surface in a way that is completely independent of the parametrization.
Along with the projectors, presented in~\cite{dassi2017}, the work in~\cite{mola2014}
describes an algorithm which allows for the hierarchical refinement of an initial blocking made of quadrilateral
cells. Across each level of refinement, the cells located in the highest curvature regions are refined, until
a prescribed accuracy is reached. Figure~\ref{fig:hull_front} shows the geometry (left side) and quadrilateral grid
generated on a planing yacht hull. As can be appreciated, the grid is consistently refined in the high curvature
regions located around the double chine line, and on the spray rails. Finally, once this quadrilateral water tight
mesh is obtained, the cells are split into triangles to obtain a water tight STL triangulation, which can be
an ideal input for numerical analysis software.    

We remark that all the presented applications exemplify the employment of the
numerical pipeline proposed in the framework of the industrial POR-FESR
projects SOPHYA ``Seakeeping Of Planing Hull YAchts'', PRELICA ``Advanced methodologies for hydro-acoustic design of naval propulsion'', and FSE HEaD
``Higher Education and Development'' programme founded by European Social Fund,
in which mathLab laboratory has been involved in the last years.

\section{Parameter space dimensionality reduction}

After all the contributions from the different teams, the number of
parameters could be too big for a reasonable optimization cycle in
terms of computational time. In other cases, even if the parameters
are not too many, there could be some of them dependent on the
others. To overcome this problem it is possible to reduce the
dimensionality of the parameter space by finding an active
subspace (AS)~\cite{constantine2015active} of the target functions. This
technique ascertains whether the output of interest can be
approximated by a function depending by linear combinations of all
the original parameters. Its application has been proven successful in
several parametrized engineering
models~\cite{grey2017active,constantine2017time,demo2018isope,tezzele2018dimension}.

Now we briefly review the process of finding active subspaces of a
scalar function $f$ representing the output of interest, and depending on
the inputs $\mupar \in \mathbb{R}^m$. Let us assume
$f: \mathbb{R}^m \rightarrow \mathbb{R}$ is a scalar function
continuous and differentiable in the support of a probability density
function $\rho: \mathbb{R}^m \rightarrow \mathbb{R}^+$. We assume $f$
with continuous and square-integrable (with respect to the measure
induced by $\rho$) derivatives. We define the active subspaces of the
pair $(f, \rho)$ as the eigenspaces of the covariance matrix
associated to the gradients $\nabla_{\mupar} f$.
The elements of this matrix, the so-called uncentered covariance matrix of the
gradients of $f$, denoted by $\mathbf{C}$, are the average products of partial
derivatives of the simulations' input/output map, i.e.:
\begin{equation*}
\label{eq:covariance}
\mathbf{C} = \mathbb{E}\, [\nabla_{\mupar} f \, \nabla_{\mupar} f
^T] = \int_{\mathbb{D}} (\nabla_{\mupar} f) ( \nabla_{\mupar} f )^T
\rho \, d \mupar ,
\end{equation*}
where $\mathbb{E}[\cdot]$ is the expected value.
The matrix $\mathbf{C}$
is symmetric and positive semidefinite, so it admits a real eigenvalue
decomposition $\mathbf{C} = \mathbf{W} \mathbf{\Lambda} \mathbf{W}^T$,
where $\mathbf{W}$ is a $m \times m$ orthogonal matrix of eigenvectors,
and $\mathbf{\Lambda}$ is the diagonal matrix of the eigenvalues,
which are non-negative, arranged in descending order.

Now we select the first $M$ eigenvectors, for some $M < m$, forming a
lower dimensional parameter subspace. We underline that, on average,
low eigenvalues suggest the corresponding vector is in the nullspace
of the covariance matrix. So we can construct an approximation of
$f$ by taking the eigenvectors corresponding to the most energetic
eigenvalues. Let us partition $\mathbf{\Lambda}$ and $\mathbf{W}$ as
follows:
\[
\mathbf{\Lambda} =   \begin{bmatrix} \mathbf{\Lambda}_1 & \\
                                     &
                                     \mathbf{\Lambda}_2\end{bmatrix},
\qquad
\mathbf{W} = \left [ \mathbf{W}_1 \quad \mathbf{W}_2 \right ],
\]
where $\mathbf{\Lambda}_1 = \text{diag}(\lambda_1, \dots, \lambda_M)$, and
$\mathbf{W}_1$ contains the first $M$ eigenvectors. The range of
$\mathbf{W}_1$ is the active subspace, while the inactive subspace is
the range of $\mathbf{W}_2$. By projecting the full parameter space
onto the active subspace we approximate the behaviour of the target
function with respect to the new reduced parameters.

The active variable $\mupar_M$, and the inactive one $\etapar$, are
obtained from the input parameters as follows: 
\begin{equation*}
\label{eq:active_var}
\mupar_M = \mathbf{W}_1^T\mupar \in \mathbb{R}^M, \qquad
\etapar = \mathbf{W}_2^T \mupar \in \mathbb{R}^{m - M} .
\end{equation*}
That means that we can express any point in the parameter space $\mupar \in
\mathbb{R}^m$ in terms of $\mupar_M$ and $\etapar$ as:
\[
\mupar = \mathbf{W}\mathbf{W}^T\mupar =
\mathbf{W}_1\mathbf{W}_1^T\mupar +
\mathbf{W}_2\mathbf{W}_2^T\mupar = \mathbf{W}_1 \mupar_M +
\mathbf{W}_2 \etapar.
\]
So we can rewrite $f$ as $f (\mupar) =  f (\mathbf{W}_1 \mupar_M + \mathbf{W}_2 \etapar)$,
and construct a surrogate quantity of interest $g$ using only the
active variable $\mupar_M$
\[
f (\mupar) \approx g (\mathbf{W}_1^T \mupar) = g(\mupar_M).
\]

Active subspaces can also be seen in the more general context of ridge
approximation~\cite{constantine2016near,pinkus2015ridge,keiper2015analysis}.

\begin{figure}[ht]
\centering
\subfloat[Original configuration.]{\includegraphics[width=0.35\textwidth]{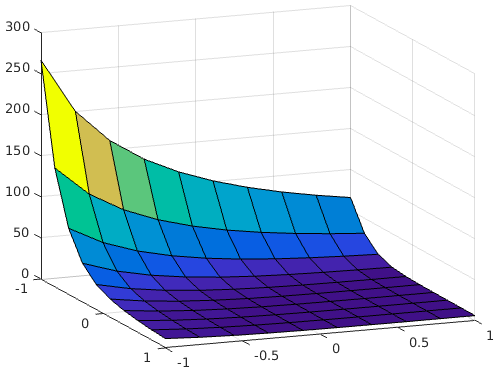}\label{subfig:as_1}}
\qquad\quad
\subfloat[The function with respect to the active variable.]{\includegraphics[width=0.35\textwidth]{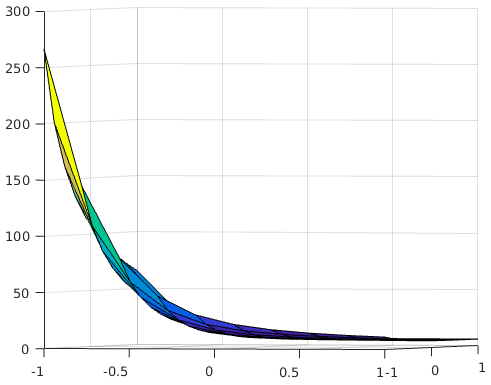}\label{subfig:as_2}}
\caption{Example of a scalar output function depending on two
  parameters. After a proper rotation of the whole domain it is
  possible to highlight the behaviour of the quantity of interest
  along the active direction.}
\label{fig:as_example}
\end{figure}

From a practical point of view, expressing a target function with
respect to its active variable means to rescale the parameter space on
the origin and then rotate it so as to unveil the low rank behaviour
of the function. In Figure~\ref{fig:as_example} an example
involving a bivariate scalar function is depicted. This approach can be viewed as a
preprocessing step in the optimization cycle that helps both in
reducing the number of the parameters and to increase the accuracy of
a further reduction of the model as proven
in~\cite{tezzele2018combined} for the computation of the pressure drop in
an occluded carotid artery using active subspace and POD-Galerkin
methods. 

Other approaches focus only on the shape parameters. To retain a
significant geometric variance while reducing the number of
geometrical parameters there exist nonlinear methods such as, among others,
Kernel PCA~\cite{scholkopf1997kernel,scholkopf1998nonlinear}, Local
PCA~\cite{kambhatla1997dimension}, and particular neural network such
as Auto Encoders and Deep Auto Encoders~\cite{hinton2006reducing}. For
a comprehensive comparison among them we refer
to~\cite{d2017nonlinear}, where the methods are demonstrated for the
design-space dimensionality reduction of a destroyer hull. For a
comparison between 13 different nonlinear techniques
see~\cite{van2009dimensionality}.

\section{Data driven model order reduction}
\label{sec:dmd}

In the \emph{big data} era, data-driven models is becoming more and more
popular in order to extract as much information as possible from all the data
acquired during the physical experiments and the simulations.
We mention also \emph{uncertainty quantification} and ROM algorithms modified
``ad hoc''~\cite{ChenQuarteroniRozza2017}.
Also in model order reduction community, several techniques have been developed
to face industrial problems with a non-intrusive approach. 

\subsection{Dynamic mode decomposition}
\label{sec:dmd}
Dynamic mode decomposition (DMD) has emerged as a powerful tool for
analyzing the dynamics of nonlinear systems, and for postprocessing
spatio-temporal data in fluid
mechanics~\cite{schmid2011application,schmid2011applications,stegeman2015proper}. It
was developed by Schmid in~\cite{schmid2010dynamic}, and it is an
equation-free algorithm, and it does not make any assumptions about
the underlying system. DMD allows to describe a non-linear
time-dependent system as linear combination of few main structures
evolving linearly in time. Many variants of the DMD were developed in
the last years like forward backward
DMD, compressed DMD~\cite{erichson2016compressed}, multiresolution
DMD~\cite{kutz2016multiresolution}, higher order
DMD~\cite{le2017higher}, and DMD with
control~\cite{proctor2016dynamic} among others. For a complete review refer
to~\cite{kutz2016dynamic,tu2014dynamic}, while for an implementation
of them we refer to the Python package called PyDMD~\cite{demo2018pydmd}. Lots of these variants arose
to solve particular industrial problems such as streaming DMD~\cite{hemati2014dynamic}
that are able to feed the classical algorithm with new real-time data
coming from sensors, and do not require storage of past
data, and they prove useful for real-time PIV or
smoke/dye visualizations. In presence of very large dataset for complex industrial model
the DMD modes are computed via randomized
methods~\cite{erichson2016randomized}. We cite also a new paradigm for
data-driven modeling that simultaneously learns the dynamics and
estimates the measurement noise at each observation that uses deep
learning and DMD for signal-noise decomposition~\cite{rudy2018deep}.

Now we present a brief overview of the standard algorithm.
Let us consider $m$ vectors, equispaced in time, representing the
state of our system, also called \textit{snapshots}: $\{\boldsymbol{x}_i\}_{i=1}^{m}$.  The idea is that
there exists a linear operator $\mathbf{A}$ that approximates the
nonlinear dynamics of $\boldsymbol{x} (t)$, i.e. $\boldsymbol{x}_{k+1}
= \mathbf{A} \boldsymbol{x}_k$. Without explicitly computing the
operator $\mathbf{A}$ we seek to approximate its eigenvectors and
eigenvalues, and we call them  DMD modes and eigenvalues. First of all
we arrange the snapshots in two matrices $\mathbf{X}$ and $\mathbf{Y}$
so as each column of the latter contains the state
vector at the next timestep of the one in the corresponding
$\mathbf{X}$ column, as follows
\begin{equation*}
\label{eq:matarranged}
\mathbf{X} =
 \begin{bmatrix}
  x_1^1   & x_2^1  & \cdots & x_{m-1}^1 \\
  x_1^2   & x_2^2  & \cdots & x_{m-1}^2 \\
  \vdots  & \vdots & \ddots & \vdots    \\
  x_1^n   & x_2^n  & \cdots & x_{m-1}^n 
 \end{bmatrix},\quad\quad
 \mathbf{Y} =
 \begin{bmatrix}
  x_2^1   & x_3^1  & \cdots & x_m^1  \\
  x_2^2   & x_3^2  & \cdots & x_m^2  \\
  \vdots  & \vdots & \ddots & \vdots \\
  x_2^n   & x_3^n  & \cdots & x_m^n 
 \end{bmatrix}.
\end{equation*}
We are looking for $\mathbf{A}$ such that $\mathbf{Y} \approx
\mathbf{A} \mathbf{X}$. The best-fit $\mathbf{A}$ matrix is given by
$\mathbf{A} = \mathbf{Y} \mathbf{X}^\dagger$, where the symbol
$^\dagger$ represents the Moore-Penrose pseudo-inverse.

The DMD algorithm projects the data onto a low-rank subspace defined by the
POD modes, that are the first $r$ left-singular vectors of the matrix
$\mathbf{X}$. We compute them via truncated singular value
decomposition as $\mathbf{X} \approx \mathbf{U}_r
\boldsymbol{\Sigma}_r \mathbf{V}^*_r$. The unitary matrix
$\mathbf{U}_r$ contains the first $r$ modes. So we can express the
reduced operator $\mathbf{\tilde{A}} \in \mathbb{C}^{r\times r}$ as
\[
\mathbf{\tilde{A}} = \mathbf{U}_r^* \mathbf{A} \mathbf{U}_r =
\mathbf{U}_r^* \mathbf{Y} \mathbf{X}^\dagger \mathbf{U}_r =
\mathbf{U}_r^* \mathbf{Y} \mathbf{V}_r \boldsymbol{\Sigma}_r^{-1} \mathbf{U}_r^* \mathbf{U}_r =
\mathbf{U}_r^* \mathbf{Y} \mathbf{V}_r \boldsymbol{\Sigma}_r^{-1},
\]
avoiding the computation of the high-dimensional
operator~$\mathbf{A}$. $\mathbf{\tilde{A}}$ defines the linear evolution of the
low-dimensional model $\boldsymbol{\tilde{x}}_{k+1} =
\mathbf{\tilde{A}} \boldsymbol{\tilde{x}}_k$, 
where $\boldsymbol{\tilde{x}}_k \in \mathbb{R}^r$ is the low-rank
approximated state. The high-dimensional state $\boldsymbol{x}_k$ can then
be easily computed as $\boldsymbol{x}_k = \mathbf{U}_r
\boldsymbol{\tilde{x}}_k$.

Exploiting the eigendecomposition of $\mathbf{\tilde{A}}$, that is
$\mathbf{\tilde{A}} \mathbf{W} = \mathbf{W} \boldsymbol{\Lambda}$, we
can reconstruct the eigenvectors and eigenvalues of the matrix
$\mathbf{A}$. The elements in $\boldsymbol{\Lambda}$ correspond to the nonzero
eigenvalues of $\mathbf{A}$, while the eigenvectors of $\mathbf{A}$
can be computed in two ways. The first one is
by projecting the low-rank approximation $\mathbf{W}$ on the
high-dimensional space: $\boldsymbol{\Phi} = \mathbf{U}_r \mathbf{W}$.
We call the eigenvectors $\boldsymbol{\Phi}$ the \textit{projected} DMD
modes. The other possibility is the so called \textit{exact} DMD modes~\cite{tu2014dynamic},
  that are the real eigenvectors of $\mathbf{A}$, and are computed as
  $\boldsymbol{\Phi} = \mathbf{Y}\mathbf{V}_r \boldsymbol{\Sigma}_r^{-1}
  \mathbf{W}$.

\begin{figure}[ht]
    \includegraphics[width=\textwidth]{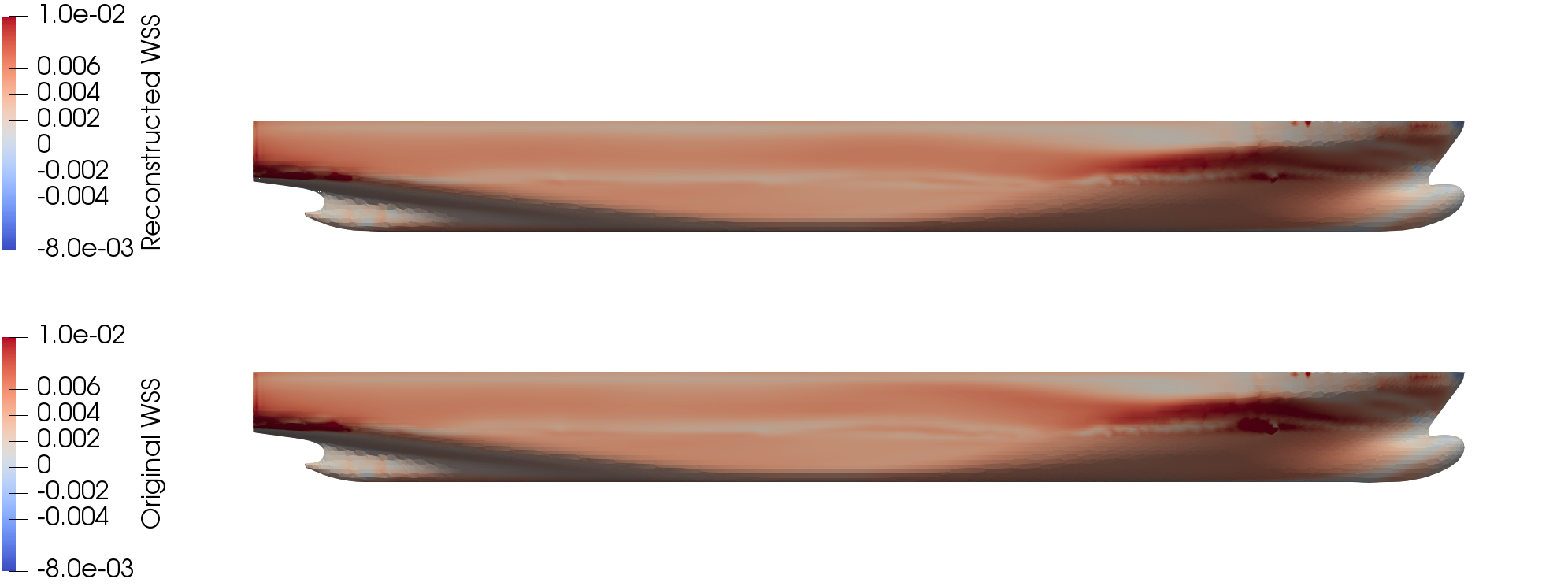}
    \caption{Example of DMD application for wall shear stress prevision. In the
    top image, we show the wall shear stress along the $x$ direction, at time
    $t = 50s$, computed using full-order solver. In the bottom, we show the
    wall shear stress along $x$ direction reconstructed at time $t = 50s$ using
    30 snapshots equispaced in the temporal window $[1, 30]$.}
\label{fig:dmd}
\end{figure}

DMD has also been successfully used to accelerate the computation of
the total drag resistance of a hull advancing in calm
water~\cite{tezzele2018model,demo2018shape,demo2018isope}. This
responded to the industrial needs of a rapid creation of the offline dataset.
We decided not only to identify the approximated dynamics of the system but
also to predict its evolution in order to achieve the regime state using only
few snapshots, as we show in the example reported in Figure~\ref{fig:dmd}.

\subsection{Proper orthogonal decomposition with interpolation}
\label{sec:podi}
Proper orthogonal decomposition with interpolation (PODI) is an equation-free model
order reduction technique providing a fast approximation of the solution of a
parametric PDE. The key idea is to approximate the solution manifold by
interpolating a finite set of high-fidelity snapshots, computed for some chosen
parameters. Since interpolation of high dimensional data can be very expensive,
    we need reduced order modelling for a real-time evaluation of
    the solution for the new parameters.

This method consists in two logical phases: in the
\textit{offline} one, the high-fidelity solutions of a finite set of deformed
configurations are computed and stored into the matrix $\mathbf{S}$ such that:
\begin{equation*}
\mathbf{S} = \begin{bmatrix} s_1 & s_2 & \dotsc & s_m \end{bmatrix},\quad s_i
\in \mathbb{R}^n\quad\text{for}\,i = 1, 2, \dotsc, m.
\end{equation*}
The basis spanning the low dimensional space is computed applying the singular
value decomposition on the snapshots matrix:
\begin{equation*}
\mathbf{S} = \mathbf{U} \mathbf{\Sigma} \mathbf{V}^*,
\end{equation*}
where $\mathbf{U} \in \mathbb{C}^{n\times m}$ refers to the matrix whose
columns are the left singular vectors --- the so called \emph{POD modes} --- of
the snapshots matrix. We project the high-fidelity solutions onto the low-rank
space, so they are represented as linear combination of the
modes and the coefficients of this combination are called modal coefficients.

In the \textit{online} phase the modal coefficients are interpolated and
finally, for any new parameter, the solution of the parametric PDE is
approximated.  This method has the great benefit of being based only on the
system output, but the accuracy of the approximation depends on the chosen
interpolation method. The algorithm has been implemented in an open source
Python package called EZyRB~\cite{demo2018ezyrb}. For deeper details about the
PODI, we recommend~\cite{salmoiraghi2018,bui2003proper,ripepi2018reduced}.

\section{Simulation-based design optimization framework}

As previously stated, a shape optimization pipeline is usually composed by three
fundamental ingredients: a deformation technique to construct the set
of admissible shape, an objective function, and an optimal
strategy to converge to the optimal shape with the lowest number of
evaluations.
Depending on the studied physical phenomena, the entire process can be very
long: many complex problems, as for example conductivity, diffusion and fluid
dynamic, are described through partial differential equations (PDEs). The
numerical solution of such equations is usually expensive from the
computational viewpoint. Moreover, in an optimization scenario, these equations
have to be solved at each iteration, making the computational cost
unaffordable for many applications, especially in the industrial sectors where
a high responsiveness is requested to reduce the time-to-market. The model
order reduction (MOR) offers the possibility to efficiently compute the
solution of parametric PDEs, drastically reducing the computational effort. We
exploited MOR techniques to design an innovative shape optimization pipeline
which fits the industrial needs, primarily in terms of efficiency, reliability
and modularity. The key idea of this optimization procedure is to collect
the solutions, or the output of interest, from the full-order model for a
finite set of parameters, then combining these solutions for a fast evaluation
of the solution for any new iteration of the optimization algorithm. 

In the first step, the deformed shapes are created from the initial geometry by
using a combination of parameters. There are many possible techniques
to choose from, as presented in Section~\ref{sec:geom}. The important
aspect is that given a set of parameters, the software is able to
generate a new deformed geometry.

The parameter space is sampled and the system
configurations so-created are evaluated using the high-fidelity numerical
method. The pipeline relies only on the system outputs, without requiring
information about the physical system, making all the procedure independent
from the high-fidelity solver. Especially in an industrial context, this
guarantees a great plus, allowing to adopt any solver --- also commercial ---
within the pipeline.  Further, the non-intrusive approach preserves
the industrial know-how and reduces the complexity in the implementation phase.

We use two different data-driven model order reduction methods to accelerate
the optimization. With the dynamic mode decomposition described in Section~\ref{sec:dmd}
we can simulate the physical
problem at hand for a shorter temporal window using the computational expensive
full-order solver and apply the DMD on the produced output to predict the
solution/output of interest at regime. The second model order reduction
technique adopted in the pipeline is the PODI, discusses in Section~\ref{sec:podi}.
Thanks to this method, we have the possibility to approximate in a real-time
context the solution of parametric PDEs, combining several pre-computed
snapshots. We adopt PODI in order to deal only with the output data of the
high-fidelity solver, thus let the pipeline be independent from the used
full-order model.  To increase the accuracy of the reduced order model, an
intrusive approach can be adopted. For an exhaustive discussion on the
intrusive model order reduction, we suggest~\cite{stabile2017,stabile2018} due
to the implementation of MOR methods in a finite volume (FV) framework, the
nowadays industrial standard for many fluid dynamics applications. For
an overview on projection-based ROMs and the effort in increasing the
Reynolds number see~\cite{HijaziAliStabileBallarinRozza2018,ballarin2015supremizer}, and
for the joint use of such methods and uncertainty quantification strategies based on
non intrusive polynomial chaos see~\cite{hijazi2018quiet}. Another
possibility is to link together the isogeometric analysis with MOR into a
complete parametric design pipeline from CAD to accurate and efficient
numerical simulation~\cite{salmoiraghi2016isogeometric,garotta2018quiet}. For a
complete discussion about ROM for parametric PDEs, we
recommend~\cite{HesthavenRozzaStamm2015}.

\begin{figure}[!t]
\centering
\includegraphics[width=.47\textwidth]{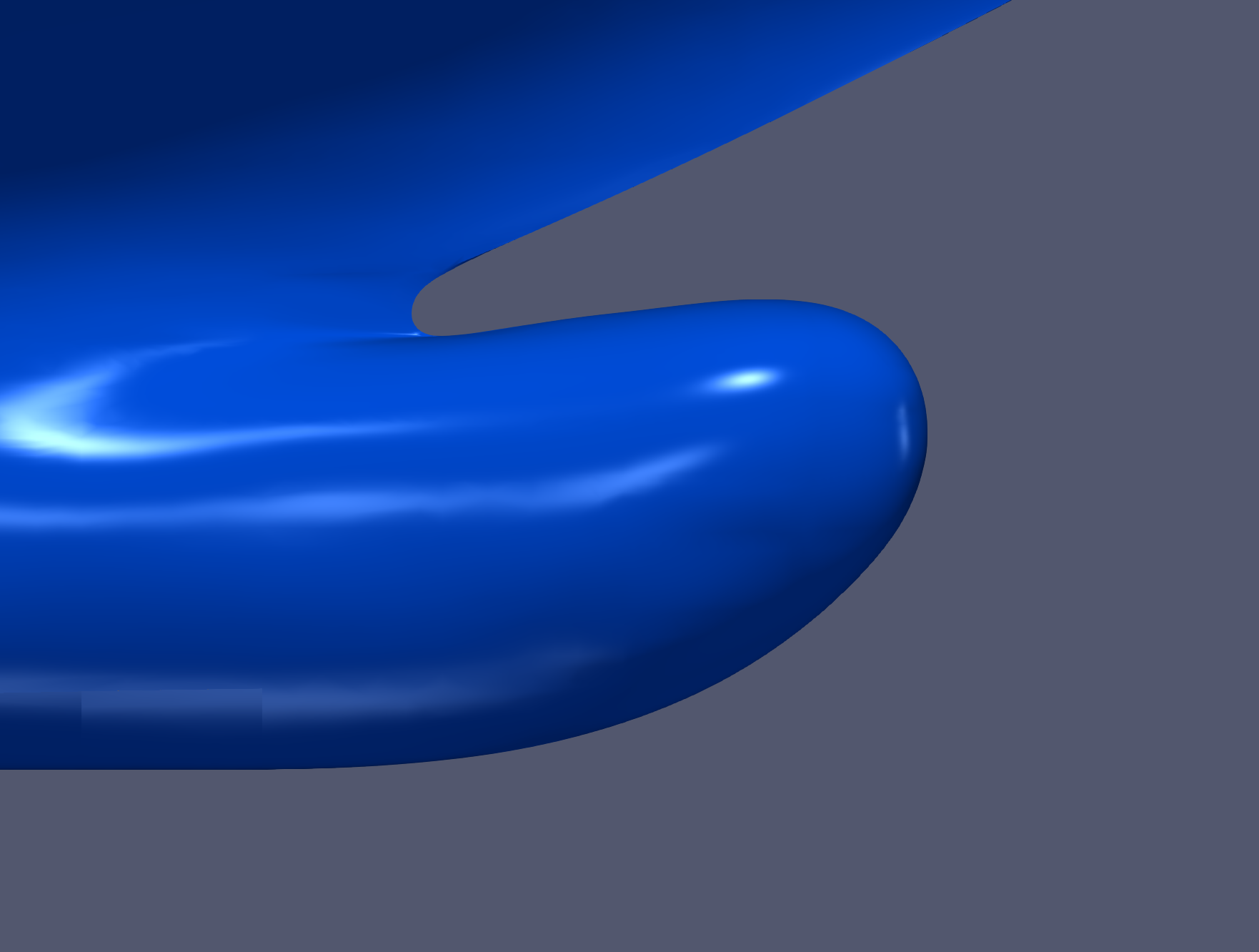} \hfill
\includegraphics[width=.47\textwidth]{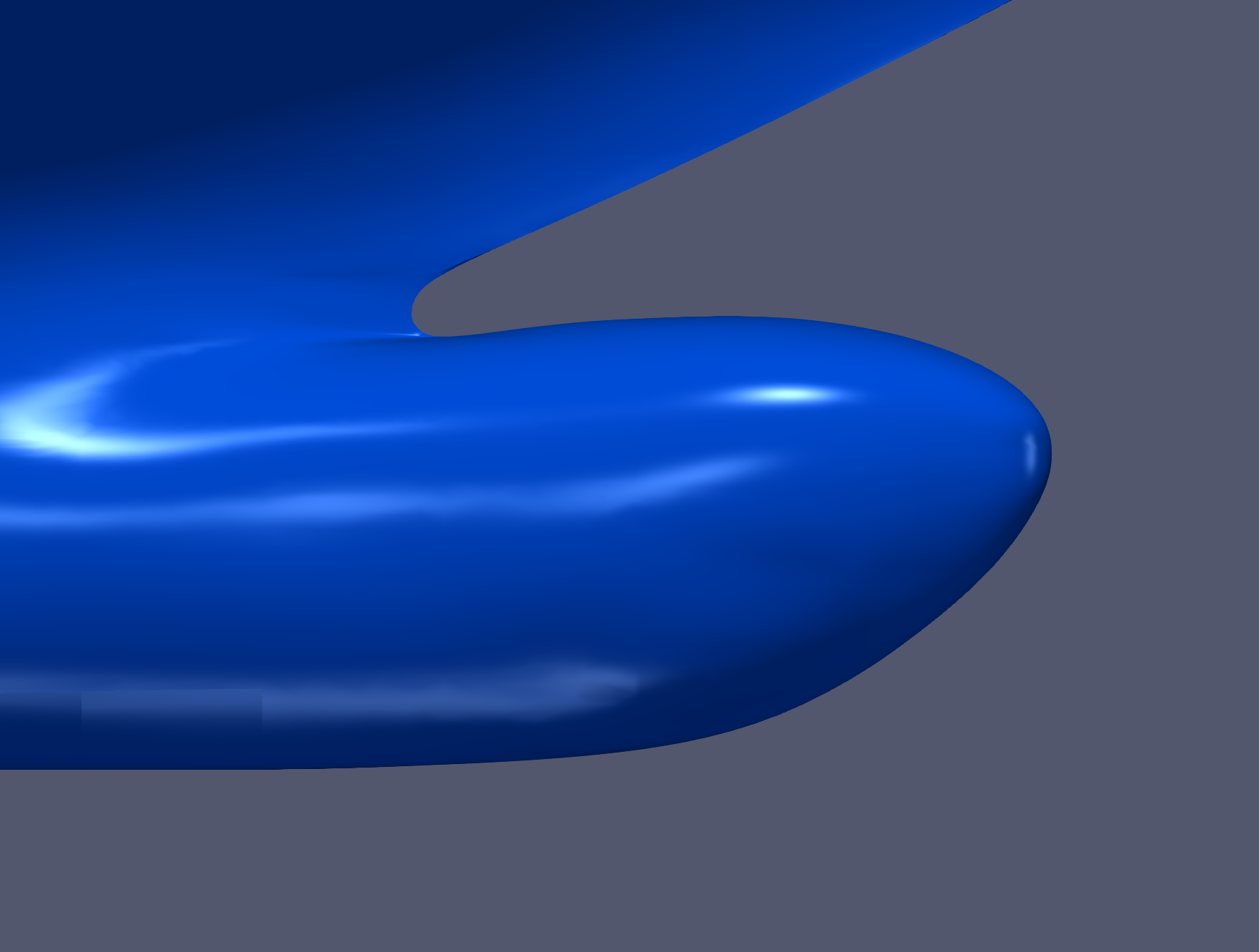}
\caption{Example of the shape optimization pipeline applied on naval
  hull: the original shape (left) and the optimized shape (right).}
\label{fig:opt}
\end{figure}

The optimization algorithm relies so on the reduced order model: since the
online phase returns the approximated solutions in a quasi real-time scenario,
the optimization algorithm lasts minutes or hours to reach the optimal shape,
also if thousands of iterations are needed. The computational cost of the
procedure is due to the creation of the solutions database.  Thanks to MOR, we
have also the possibility to run and tests many different optimization
algorithms, avoiding any further high-fidelity simulations. Moreover, the solutions database can be
enriched to increase the accuracy of the reduced order model. Examples of
optimization procedure involving MOR techniques applied into naval and
aerodynamics fields are
respectively~\cite{demo2018shape,scardigli2019enabling}. Figure~\ref{fig:opt}
shows the results of the application of the shape optimization system
on the bulbous bow of a cruise ship. This achievement has been
developed in the framework of a regional European Social Fund project
from Regione Friuli Venezia-Giulia: HEaD in collaboration with Fincantieri -
Cantieri Navali Italiani S.p.A..

\section{Conclusions and perspectives}
\label{sec:the_end}
Industrial computational needs are every day more and more demanding in terms
of computational time, reliability, error certification,
data-assimilation, robustness, and easiness of use. In this work we
presented several model order reduction and shape parameterization
techniques to solve industrial and applied mathematics problems.

\begin{figure}[h!]
\centering
\includegraphics[width=1.\textwidth]{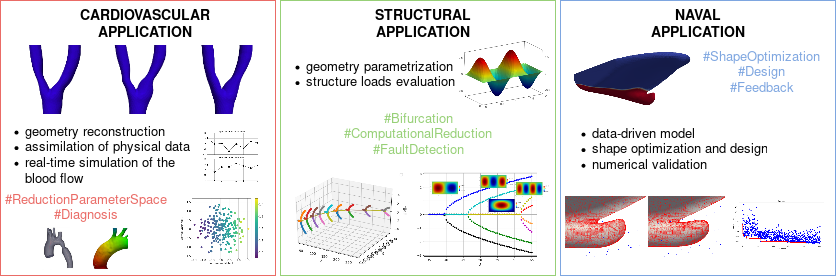}
\caption{Possible extension of the presented pipeline with different
goal and application fields. From cardiovascular problems like
real-time blood flow simulation, to structure load analysis and
identification, as well as naval applications.}
\label{fig:global_scheme}
\end{figure}

More has to be made to integrate real-time data-assimilation, machine
learning and prediction, but we are moving along this horizon and
MOR will play a crucial role to tackle many complexities arising from
complex industrial artifacts management. A step in this direction is
the planned webserver ARGOS~\cite{argos}, developed by mathLab group at SISSA
that will
make possible the exploitation of reduced order models to a vast
category of people working in design, structural, and CFD
teams. Through specific web applications the user will be able to
solve many industrial and biomedical problems without the need of
being an expert in numerical analysis and scientific computing. Figure~\ref{fig:global_scheme}
depicts some of the possible applications that are currently being developed.

We also cite the Artificial Student ``Artie''~\cite{patera2018project}
that accepts problem statements posed in natural language, and solves
numerically some PDEs problems, that will help both students, the
scientific staff, and engineers in general. Moreover we want to
highlight the effort of the Italian government in the technology
transfer thanks to the institution of several competence centers
connecting research facilities, university, and companies in the
framework of Industry 4.0. Similar initiatives are undergoing in many other
European countries (France, Germany, UK, \dots).

\section*{Acknowledgements}
We acknowledge the scientific collaboration within SISSA
mathLab group and the support provided by Dr Francesco Ballarin and Mr
Federico Pichi for casting future extensions.

This work was partially performed in the context of the project SOPHYA -
``Seakeeping Of Planing Hull YAchts'' and of the project PRELICA - ``Advanced methodologies for hydro-acoustic design of naval propulsion'', both supported by Regione
FVG, POR-FESR 2014-2020, Piano Operativo Regionale Fondo Europeo per
lo Sviluppo Regionale, partially funded by the project HEaD, ``Higher
Education and Development'' in collaboration with Fincantieri, supported by Regione FVG, European
Social Fund FSE 2014-2020, and partially supported by European Union Funding for
Research and Innovation --- Horizon 2020 Program --- in the framework
of European Research Council Executive Agency: H2020 ERC CoG 2015
AROMA-CFD project 681447 ``Advanced Reduced Order Methods with
Applications in Computational Fluid Dynamics'' P.I. Gianluigi
Rozza.



\end{document}